\documentclass[a4paper,leqno,10pt,twoside]{amsart}
\usepackage{amsmath,amsfonts,amssymb,amsthm,amscd}
\usepackage[utf8]{inputenc}
\usepackage[english]{babel}
\usepackage[colorlinks=true,citecolor=blue, urlcolor=blue, linkcolor=blue,pagebackref]{hyperref}
\usepackage[top=1.1in,bottom=1.1in,left=1.in,right=1.in]{geometry} 
\usepackage{graphicx}
\usepackage{wrapfig}
\usepackage{paralist}
\usepackage{tabto}
\usepackage{standalone}
\usepackage{xfrac}
\usepackage{float}
\usepackage{tikz-cd}
\usepackage{mathtools}

\usepackage{pgf,tikz}
\usetikzlibrary{arrows,chains,
 decorations.pathreplacing,
shapes,%
}

\usepackage[all]{xy} 

\usepackage{enumerate}
\usepackage{xcolor}
\usepackage{aliascnt}

\theoremstyle{plain}
\newtheorem{thm}{Theorem}[section]

\newtheorem{thmIntr}{Theorem}

\newaliascnt{propIntr}{thmIntr}

\newaliascnt{corIntr}{thmIntr}

\newaliascnt{QU}{thm}

\aliascntresetthe{QU}

\newaliascnt{lem}{thm}
\newtheorem{lem}[lem]{Lemma}
\aliascntresetthe{lem}

\newaliascnt{cor}{thm}
\newtheorem{cor}[cor]{Corollary}
\aliascntresetthe{cor}

\newaliascnt{prop}{thm}
\newtheorem{prop}[prop]{Proposition}
\aliascntresetthe{prop}

\theoremstyle{definition}

\newaliascnt{rem}{thm}
\newtheorem{rem}[rem]{Remark}
\aliascntresetthe{rem}

\newaliascnt{defn}{thm}
\newtheorem{defn}[defn]{Definition}
\aliascntresetthe{defn}

\newaliascnt{setup}{thm}
\newtheorem{setup}[setup]{Set-up}
\aliascntresetthe{setup}

\newaliascnt{ex}{thm}

\aliascntresetthe{ex}

\numberwithin{equation}{section}

\DeclareMathOperator{\coker}{coker}

\DeclarePairedDelimiter\ceil{\lceil}{\rceil}
\DeclarePairedDelimiter\floor{\lfloor}{\rfloor}

\definecolor{applegreen}{rgb}{0.55, 0.71, 0.0}

\title[The polarized degree of irrationality of $K3$ surfaces]{The polarized degree of irrationality of $K3$ surfaces}
\author[F.~Moretti]{Federico Moretti}

\address{F.M.: Institut für Mathematik, Humboldt-Universität zu Berlin, Unter den Linden 6, 10099 Berlin, Germany}
\curraddr{Stony Brook University,  100 Nicolls Road, Stony Brook, NY 11794}
\email{federico.moretti@stonybrook.edu}

\begin{document}

\keywords{Degree of irrationality, $K3$ surfaces, kernel bundle, abelian varieties}
\subjclass[2020]{14E05, 14F08, 14J28}

\maketitle

\setcounter{tocdepth}{1}

	\begin{abstract}
		Given a polarized variety $(X,L)$,  we construct and study projections of low degree $ X\dashrightarrow \mathbb{P}(H^0(L^\vee)) \dashrightarrow  \mathbb P ^n $ using the associated kernel bundles.  As an application, we can show that the degree of irrationality of a very general $(1,6)-$polarized abelian surface, as well as that of a very general $K3$ surface of genus $6$ is $3$. We also give new upper bounds on the degree of irrationality of $K3$ surfaces of any genus. We study the family of projections of minimal degree of a very general $K3$ surface of genus $4,5,6$. As a different application of our construction, we exhibit new generically finite rational maps of low degree from some hyper-K\"ahler varieties and abelian varieties to projective spaces.
	\end{abstract}

	\section*{Introduction}

	Given a quasi-projective variety $X$ over an algebraically closed field, a natural invariant measuring how far $X$ is from being rational is its \emph{degree of irrationality}
\[
\mathrm{irr}(X)=\min\Bigl\{\deg(\varphi)\ \Big|\ \varphi:X \dashrightarrow \mathbb{P}^{\dim(X)} \text{ is a generically finite rational map}\Bigr\}.
\]

The degree of irrationality of a curve is called its \emph{gonality} and is a classical invariant. In contrast, when $\dim(X)\ge 2$ very little is known in general, and the problem has only recently received considerable attention. The invariant $\mathrm{irr}(X)$ was introduced in \cite{foundation}. This invariant has been studied in a number of cases, including hypersurfaces \cite{4,5,7,St,22}, abelian varieties \cite{AlzatiPirola,AS,av,absurf,35,21,yoshi}, hyper-K\"ahler varieties \cite{VS}, and other classes of examples \cite{ago,fano}. For a recent survey the reader may consult \cite{MR5049528}. Giving sharp lower bounds is, in general, a difficult problem. In this paper, among other things, we compute the degree of irrationality of a very general $K3$ surface of genus $6$ and of a very general $(1,6)$-polarized abelian surface. The latter was one of the remaining open cases for very general $(1,d)$-polarized abelian surfaces in view of \cite{AS,absurf}. The only other case in which the degree of irrationality of a very general $(1,d)$-polarized abelian surface was known to be $3$ is $d=2$, as proven in \cite{35}.

More generally, we introduce and study a natural invariant of polarized projective varieties closely related to the degree of irrationality. Let $(X,L)$ be a polarized projective variety of dimension $n\ge 2$, where $L$ is an ample line bundle. For any subspace $V\subset H^0(L)$ of dimension $n+1$ we obtain a rational map
\[
\varphi_V:X \dashrightarrow \mathbb{P}^n.
\]
We define
\[
\mathrm{irr}_L(X)=\inf\Bigl\{\deg(\varphi_V)\ \Big|\ V\in \mathrm{Gr}(n+1,H^0(L)),\ V \text{ generates $L$ in codimension }1 \text{ and $\varphi_V$ is generically finite}\Bigr\}.
\]
Equivalently, given the map $X\dashrightarrow \mathbb P H^0(L)^\vee$, we aim to study the projections $\mathbb P H^0(L)^\vee \dashrightarrow \mathbb P^n$ such that the induced map $X\dashrightarrow \mathbb P^n$ is generically finite and has minimal degree.

\textbf{The curve case.}
Let $C$ be a general curve of genus $g$, and let $L\in \mathrm{Pic}^e(C)$ be a general line bundle of degree $e\ge g+1$. In this case it is elementary to see that
\[
\mathrm{irr}_L(C)=\min\Bigl\{\, d \ \Big|\ L\in W^1_d(C)+W^0_{e-d}(C)\Bigr\},
\]
where $W^r_d(C)\subset \mathrm{Pic}^d(C)$ denotes the Brill--Noether locus of line bundles of degree $d$ carrying at least $r+1$ linearly independent global sections and $W^1_d(C)+W^0_{e-d}(C)$ is the image of the summation map
\[
W^1_d(C)\times W^0_{e-d}(C)\longrightarrow \mathrm{Pic}^e(C),\qquad (A,B)\longmapsto A\otimes B.
\]
Since $L$ is general and the summation map is of maximal rank, the condition $L\in W^1_d(C)+W^0_{e-d}(C)$ is equivalent to requiring that the locus $W^1_d(C)+W^0_{e-d}(C)$ has dimension $g$ (i.e.\ dominates $\mathrm{Pic}^e(C)$). Thus $\mathrm{irr}_L(C)$ is the minimal $d$ for which
\[
\dim W^1_d(C)+\dim W^0_{e-d}(C)\ \ge\ g.
\]
By Brill--Noether theory, for a general curve one has
\[
\dim W^1_d(C)=\mathrm{min}\{\rho(g,1,d),g\}=\mathrm{min}\{2d-g-2,g\},
\qquad
\dim W^0_{e-d}(C)=\mathrm{min}\{e-d,g\}.
\] Plugging these expressions into the inequality above yields $\mathrm{irr}_L(C)=\mathrm{max}\{2g+2-d,\lfloor{\frac {g+3} 2}\rfloor\}$.

For a surface $S$ and a non-trivial indecomposable class $L\in \mathrm{Pic}(S)$ (i.e.\ $L$ cannot be written as the tensor product of two non-trivial effective line bundles), for any $V \in \mathrm{Gr}(3,H^0(L))$  the rational map $\varphi_V$ is generically finite and $V$ generates $L$ in codimension $1$. For $K3$ surfaces we obtain the following.

\begin{thmIntr}\label{k3}
Let $(S,L)$ be a $K3$ surface of genus $g=2+2n(n+1)+k$. Then
\[
\mathrm{irr}_{L}(S)\le  2+n +\ceil*{\frac k 2}- \floor*{ \frac k 4 }.
\]
Moreover, if $g=3,4,6$ and $S$ is very general, then $\mathrm{irr}_L(S)=3$. If $g=5$ and $S$ is very general, then $\mathrm{irr}_L(S)=4$.
\end{thmIntr}

In the paper we present the bound for $K3$ surfaces of Picard rank $1$, and the corresponding statement for arbitrary $K3$ surfaces follows from general specialization results; see \cite[Corollary~C]{7}. The bounds we obtain improve those in \cite{St}. In low genus they are optimal up to genus $8$ (in \cite{morettirojas} it is proven that $\mathrm{irr}_L(S)=4$ for $g=7,8,9,11$), and asymptotically the improvement over \cite{St} is by a factor between $3$ and $6$, depending on the genus. Let us also point out that after the appearance of the first version of this paper the bounds were further improved: in \cite{morettirojas} the bounds were sharpened for $9\le g\le 14$ by combining the method of the present paper with derived category tools (Fourier--Mukai transforms and Bridgeland stability conditions), while in \cite{Oyanedel} our approach was optimized to yield better asymptotic bounds for $k\gg 0$ in \autoref{k3}.

Our technique relies on an elementary yet powerful observation, which we briefly outline below. It can be seen as a natural variation of the Lazarsfeld--Mukai bundles introduced to study linear series on curves on $K3$ surfaces; see \cite{LAZ}. Let $X$ be a variety of dimension $n$ and let $|V|\subset |L|$ be a linear system inducing a rational map $\varphi_{V}:X \dashrightarrow \mathbb P^{n}$. One may then construct the associated kernel reflexive sheaf, fitting into an exact sequence
\[
\begin{tikzcd}
0 \arrow{r} & E^\vee \arrow{r} & V \otimes \mathcal{O}_X \arrow{r} & L .
\end{tikzcd}
\]
Via the identification $\mathrm{Gr}(n,V)=\mathbb{P}V^\vee \subset \mathbb{P}(H^0(E))$, the fibers of $\varphi_V$ can be described as zero loci of sections of $E$. More precisely, the fiber above $[s]\in \mathbb PV^\vee=\mathbb P^n$ is $Z(s)$ outside the base locus of the map (see \autoref{wow}). This provides a quick way to compute the degree of $\varphi_V$.

The main idea of this paper is that, under suitable hypotheses on $E$, this construction can be reversed. Namely, if one is given a sufficiently general vector bundle $E$ of rank $n$ together with at least $n+1$ linearly independent global sections generating $E$ in codimension $1$, then one obtains a generically finite rational map of degree $\le c_{n}(E)$. A simple corollary, proved along these lines in \cite{Morettirojas2}, is that if $X$ carries a globally generated vector bundle $E$ with $c_n(E)>0$, then $\mathrm{irr}(X)\le c_n(E)$.

Moreover, if $E$ has many global sections, one can impose base points and make the degree drop. The basic example to keep in mind is a general $K3$ surface $S$ of genus $6$. In this case there is a stable vector bundle $E$ of rank $2$ on $S$ with $h^0(E)=5$ and $c_2(E)=4$. Then, for any point $P\in S$, the map defined by
\[
V=H^0(E \otimes \mathcal I_P)^\vee \subset H^0(L)
\]
has degree $c_2(E)-1=3$; see also \autoref{genus6}. More generally, \autoref{k3} is a direct application of \autoref{wow}. We now outline some further applications.

For a surface $S$ with $\mathrm{NS}(S)=\mathbb Z \cdot L$, this technique can be used to give a scheme structure to the locus of special maps; see \autoref{scheme}. We endow the locus
\[
W^r_d(S,L)=\Bigl\{V \in \mathrm{Gr}(r+1,H^0(L)) \ \Big| \ \deg(\varphi_V)\le d\Bigr\}
\]
with a natural structure of algebraic variety. Some explicit computations are as follows.

\begin{thmIntr}\label{456}
Let $(S,L)$ be a $K3$ surface of genus $g$ with $\mathrm{Pic}(S)=\mathbb Z \cdot L$. The following holds:
\begin{enumerate}
\item if $g=4$, then $W^2_3(S,L)=\mathbb P^3$;
\item if $g=5$, then $W^2_3(S,L)=\emptyset$ and $W^2_4(S,L)$ has two components, one birational to $S$ and one birational to a $\mathbb P^3$-bundle over the moduli space of stable rank $2$ vector bundles with $c_2=4$;
\item if $g=6$, then $W^2_3(S,L)=S$.
\end{enumerate}
\end{thmIntr}

A more geometric interpretation of \autoref{456} is the following. If $S\subset \mathbb P^4$ is a very general $K3$ surface of degree $6$, then the locus of trisecant lines\footnote{A linear space is $k$-secant to $S$ if its scheme theoretic intersection with $S$ is of length $k$.} to $S$ is isomorphic to $\mathbb P^3$. If $S\subset \mathbb P^5$ is a very general $K3$ surface of degree $8$, then there are no $4$-secant planes tangent to $S$. If $S\subset \mathbb P^6$ is a very general $K3$ surface of degree $10$, then for any $p\in S$ there exists a unique $6$-secant codimension-$3$ linear subspace tangent to $S$ at $p$. A similar theorem for $K3$ surfaces up to genus $14$ is developed in \cite[Theorem~A and Theorem~B]{morettirojas}, suggesting interesting patterns for the dimension of the loci $W^2_d(S,L)$; see \cite[Section~8]{morettirojas}. Other kinds of polarized surfaces are also interesting from this perspective: for example, the case of $(S,3\Theta)$ for a very general principally polarized abelian surface $(S,\Theta)$ should be accessible with the methods of the present paper.
	
	\, We give a few elementary applications to other kinds of varieties. For a surface $S$ we denote by $S^{[n]}$ the Hilbert scheme of zero-dimensional subschemes of length $n$.

\begin{thmIntr}\label{higherdim}
The following bounds hold:
\begin{enumerate}
\item if $(S,L)$ is a very general $(1,6)$-polarized abelian surface, then $\mathrm{irr}(S)=\mathrm{irr}_L(S)=3$;
\item if $(S,L)$ is a very general $K3$ surface of genus $6$, then $\mathrm{irr}(S^{[2]})\le 6$;
\item if $(A,L)$ is a very general abelian threefold of type $(1,3,12)$ or $(1,6,6)$, then $\mathrm{irr}_L(A)\le 8$;
\item if $(S,L)$ is a very general $K3$ surface of genus $10$, then $\mathrm{irr}(S^{[3]})\le 20$.
\end{enumerate}
\end{thmIntr}

\, Let us remark that the case of $(1,6)$-polarized abelian surfaces settles one of the remaining exceptional cases for abelian surfaces. More precisely, the degree of irrationality of a very general $(1,d)$-polarized abelian surface is either $3$ or $4$ for every positive integer $d$. It is known to be $4$ when $d\neq 1,2,3,6$ (see \cite{AS,absurf}), and it is known to be $3$ when $d=2$ (see \cite{35}). The first item in \autoref{higherdim} addresses the case $d=6$. In a somewhat similar spirit, Mason \cite{mason} completed the classification of the degree of irrationality of hyperelliptic surfaces initiated by Yoshihara \cite{yoshihara1996degrees}.

For very general abelian threefolds, the best lower bound currently known is $5$ (see \cite{av}). In \cite{Morettirojas2}, Picard bundles on Jacobians are studied, yielding the upper bound $2^g$ for the degree of irrationality of genus $g$ Jacobians.

For hyper-K\"ahler varieties of dimension $2n$ Voisin shows in \cite{VS} that the degree of irrationality is at least $n+1$.

\textbf{Structure of the paper.} In \autoref{Section1} we explain how to construct and compute the degree of a generically finite rational map starting from the kernel sheaf. In \autoref{Section2} the technique is applied to study projections of the lowest possible degree for a very general $K3$ surface $S\subset \mathbb P^g$ of genus $4,5,6$ (see \autoref{456}). In \autoref{section3} we study the case of $(1,6)$-polarized abelian surfaces. In \autoref{section4} we obtain the asymptotic bound of \autoref{k3}. We prove \autoref{higherdim} in \autoref{section5}. Finally, in \autoref{section6} we prove that the loci $W^r_d(S,L)$ are algebraic varieties (see \autoref{scheme}).

\, Each section presents an independent application of the method introduced in \autoref{Section1}, and the sections are mutually independent. The contents of \autoref{Section1} and \autoref{section3} hold over any algebraically closed field. The contents of the other sections hold over $\mathbb C$.

\textbf{Acknowledgements.} I would like to express my sincere gratitude to my advisor, Gavril Farkas, for his support, as well as for valuable discussions and his insights on the topic. This project was first inspired by the conference \emph{Algebraic Geometry in Roma 3} and a subsequent conversation with Pietro Pirola. I am grateful to the organizers of the conference and to Pietro Pirola for the stimulating exchange.

I have benefited from many insightful discussions, particularly within the algebraic geometry group at HU Berlin. I am thankful to all its members for their input. In particular, I would like to thank Andr\'es Rojas for his help with abelian varieties, as well as for numerous suggestions that significantly improved the presentation.

Further thanks go to Alessio Cela, Nathan Chen, Paolo Grossi, Andreas Leopold Knutsen, Rob Lazarsfeld, Olivier Martin, Pedro Montero, Joaquim Ribeiro and Angel Rios Ortiz for valuable discussions and helpful comments.

The author wishes to thank the referee for pointing out several inaccuracies and for providing comments that significantly improved the presentation.

This research was supported by the ERC Advanced Grant SYZYGY. It has received funding from the European Research Council (ERC) under the European Union's Horizon 2020 research and innovation programme (grant agreement No.~834172).
\section{Constructing rational maps via the kernel bundle.}
\label{Section1}
We begin with some notation and definitions.

\begin{setup}\label{setup:cn}
We denote by
\begin{itemize}
    \item $X$ a smooth projective variety of dimension $n$;
    \item $E$ a reflexive sheaf of rank $n$ on $X$;
    \item $V^\vee$ an $(n+1)$-dimensional subspace of $H^0(E)$, i.e.\ $V^\vee\in \mathrm{Gr}(n+1,H^0(E))$;
    \item $\mathrm{Bl}(V^\vee)\subset X$ the union of the locus where the sections of $V^\vee$ do not generate $E$ and the locus where $E$ is not locally free.\footnote{We will sometimes write $\mathrm{Bl}$ instead of $\mathrm{Bl}(V^\vee)$ when the dependence on $V^\vee$ is clear from the context.}
\end{itemize}
For a general section $s\in V^\vee$, let $Z_{\mathrm{cycle}}(s)$ denote the zero-cycle associated with the zero-dimensional part of $Z(s)$.

\smallskip
We define $A_s$ to be the maximal effective subcycle of $Z_{\mathrm{cycle}}(s)$ supported on $X\setminus \mathrm{Bl}(V^\vee)$. We then set
\[
c_n(V^\vee):=\deg(A_s).
\]
Note that, since $V^\vee$ generates $E$ away from $\mathrm{Bl}(V^\vee)$, a general section $s\in V^\vee$ has a zero-dimensional vanishing locus on $X\setminus \mathrm{Bl}(V^\vee)$, so $c_n(V^\vee)$ is well defined.

\smallskip
If, moreover, a general section $s\in V^\vee$ vanishes along a zero-dimensional subscheme, we define $bl(V^\vee)$ as the degree of the complementary part supported on $\mathrm{Bl}(V^\vee)$. Namely, writing
\[
Z_{\mathrm{cycle}}(s)=A_s + A_s',
\]
where $A_s'$ is supported on $\mathrm{Bl}(V^\vee)$, and set
\[
bl(V^\vee):=\deg(A_s').
\]
\end{setup}

\begin{rem}
In general, the support of $A_s'$ is a proper subset of $\mathrm{Bl}(V^\vee)$.
\end{rem}

We will be interested in \emph{pairs} $(E,V^\vee)$ that satisfy some special properties. For the purposes of this paper we will refer to them as \emph{good pairs}.

\begin{defn}\label{def:goodpair}
A pair $(E,V^\vee)$, with $E$ a reflexive sheaf of rank $n$ and $V^\vee\in \mathrm{Gr}(n+1,H^0(E))/\mathrm{Aut}(E)$, is a \emph{good pair} if $V^\vee$ generates $E$ in codimension $1$ and $H^0(E^\vee)=0$.
\end{defn}
Given a good pair $(E,V^\vee)$, the kernel of the evaluation $\mathrm{ev}:V^\vee\otimes \mathcal O_X\to E$ is $L^\vee$, where $L=\mathrm{det}(E)$. Indeed, since the evaluation is surjective in codimension $1$, we get $c_1(\mathrm{Im}(ev))=c_1(E)$ and $\mathrm{Im}(ev)^\vee=E^\vee$ (recall that $E$ is reflexive). Dualizing, we get an exact sequence \[
\begin{tikzcd}
    0 \arrow{r} & E^\vee \arrow{r} & V \otimes \mathcal O_X \arrow{r} & L.
\end{tikzcd}
\]
In particular, to any good pair $(E,V^\vee)$, we can associate the rational map induced by $( L, V)$ (notice that since $h^0(E^\vee)=0$ we have $V \subset H^0(L),$ therefore, the map is non-degenerate). We will denote the induced rational map by $\varphi_V:X\dashrightarrow \mathbb PV^\vee.$ It is easy to see that the set-theoretic support of the base locus of the rational map coincides with $\mathrm{Bl}(V^\vee)$. For instance, one can observe that $V$ coincides with the image of $\bigwedge^nV^\vee\to H^0(L)$ induced by the determinant map on the global sections of $E$. In fact, the locus where $V$ does not generate $L$ coincides with the union of the locus where $V^\vee$ does not generate $E$ and the locus where $E$ is not locally free.

\,On the other hand, as mentioned in the introduction, if one starts with a non-degenerate rational map $X\dashrightarrow \mathbb P^n$ (hence, a line bundle $L\in \mathrm{Pic}(X)$ together with $V\in \mathrm{Gr}(n+1,H^0(L))$ generating $L$ in codimension $1$), one may construct the kernel sheaf $E^\vee=\mathrm{ker}(V\otimes \mathcal O_X\to L)$ together with $V^\vee\subset H^0(E)$. We have the following.
\begin{lem}
    Given a non-degenerate rational map $X\dashrightarrow \mathbb P^n$ the pair $(E,V^\vee)$ is a good pair.
    \label{lemgoodpair}
\end{lem}

\begin{proof}
On a smooth variety any rational map is well defined in codimension $1$, hence the sequence \[
\begin{tikzcd}
    0 \arrow{r} & E^\vee \arrow{r} & V \otimes \mathcal O_X \arrow{r} & L,
\end{tikzcd}
\]
is exact in codimension $1$. Dualizing, we get that $V^\vee$ generates $E$ in codimension $1$. Moreover, $E$ is the kernel of a morphism between vector bundles, so it is reflexive.  The morphism $V\to H^0(L)$ is injective, this implies $H^0(E^\vee)=0$. We showed that $(E,V^\vee)$ is a good pair.
\end{proof}

We can now prove the key proposition that will be used throughout the paper: there is a bijective correspondence between rational maps of degree $d$  and good pairs $(E,V^\vee)$ with $c_n(V^\vee)=d$. The invariant $c_n(V^\vee)$  will make the computation of the degree of certain rational maps  much easier.  
\newpage
	\begin{prop}
		Let $X$ be a smooth projective variety of dimension $n$. There is a \emph{one-to-one} correspondence\[
        \big \{ \textrm{non-degenerate rational maps } X\dashrightarrow \mathbb P^n \textrm{ of degree $d$} \big \} \leftrightarrow \big \{ \textrm{good pairs } (E,V^\vee) \textrm{ with } c_n(V^\vee)=d \big \},
        \]
        \[
       \quad \quad \quad  \quad  \quad  \quad \quad \, \,  \, \, \, \, (L,V)\longrightarrow (E,V^\vee).
        \]
       Moreover, the fibers of the rational map $\varphi_V:X\dashrightarrow \mathbb PV^\vee$ are computed by $\varphi_V^{-1}[s]\cap (X-\mathrm{Bl})=Z(s)\cap (X-\mathrm{Bl})$, where $s$ can be seen both  as a point of the target projective space (left hand side of the equality) and as a section of the vector bundle $E$ (right hand side of the equality). In particular, if $E$ is a vector bundle and a general element of $V^\vee$ vanishes along a zero-dimensional subscheme, then
\[
c_n(V^\vee)=\deg(\varphi_V)= c_n(E)-bl(V^\vee).
\]
\label{wow}
\end{prop}

\begin{proof}
The map $(E,V^\vee)\mapsto (L,V)$ is clearly a bijection between good pairs and non-degenerate rational maps $X\dashrightarrow \mathbb P^n$. Indeed, given a rational map one can reconstruct the kernel sheaf
\[
E^\vee=\ker\bigl(V\otimes \mathcal O_X\to L\bigr)
\]
together with $V^\vee\subset H^0(E)$ (uniquely up to the action of $\mathrm{Aut}(E)$ on $H^0(E)$), and \autoref{lemgoodpair} shows that $(E,V^\vee)$ is a good pair.

It now suffices to show the second part of the proposition, as it implies that rational maps of degree $d$ correspond to good pairs $(E,V^\vee)$ with $c_n(V^\vee)=d$. We abuse notation and denote by $\varphi_V$ the induced morphism
\[
\varphi_V:(X\setminus \mathrm{Bl})\to \mathbb P V^\vee .
\]

The Euler sequence
\[
\begin{tikzcd}
0\arrow{r} & \mathcal O_{\mathbb P V^\vee}(-1) \arrow{r }& V^\vee \otimes \mathcal O_{\mathbb P V^\vee}\arrow{r} &  T_{\mathbb P V^\vee}(-1)\arrow{r} & 0
\end{tikzcd}
\]
pulls back to
\[
\begin{tikzcd}
0 \arrow{r} & L^\vee_{|_{X\setminus \mathrm{Bl}}}\arrow{r} & V^\vee \otimes \mathcal O_{X\setminus \mathrm{Bl}}\arrow{r} & E_{|_{X\setminus \mathrm{Bl}}} \arrow{r} & 0
\end{tikzcd}
\]
on $X\setminus \mathrm{Bl}$. In particular,
\[
E_{|_{X\setminus \mathrm{Bl}}}\cong \varphi_V^*T_{\mathbb P V^\vee}(-1),
\]
and pullback induces an identification
\[
\varphi_V^*:H^0\!\left(T_{\mathbb P V^\vee}(-1)\right)\xrightarrow{\ \sim\ } V^\vee \subset H^0(E).
\]
In particular, a general section of $T_{\mathbb P V^\vee}(-1)$ vanishes at the corresponding point $[s]\in \mathbb P V^\vee$ (and at no other point, since $c_n(T_{\mathbb P V^\vee}(-1))=1$). Hence the fiber over $[s]$ of the map $\varphi_V:X\setminus \mathrm{Bl}\to \mathbb P V^\vee$ is $Z(\varphi_V^*s)$.

The degree statement follows because $c_n(V^\vee)$ is, by definition, the degree of $Z(s)\cap (X\setminus\mathrm{Bl)}$ for $s\in V^\vee$ general.

Finally, the equality $c_n(V^\vee)=c_n(E)-bl(V^\vee)$ follows from the definition of $c_n(V^\vee)$ together with the standard fact that if a section of a rank-$n$ vector bundle $E$ vanishes along a zero-dimensional subscheme, then its length is exactly $c_n(E)$. Indeed, for $s\in V^\vee$ general we can write the associated zero-cycle as
\[
Z_{\mathrm{cycle}}(s)=A_s+A_s',
\]
where $A_s$ is supported on $X\setminus \mathrm{Bl}(V^\vee)$ and $A_s'$ is supported on $\mathrm{Bl}(V^\vee)$. By definition, $\deg(A_s)=c_n(V^\vee)$ and $\deg(A_s')=bl(V^\vee)$, while $\deg(Z_{\mathrm{cycle}}(s))=c_n(E)$. Therefore
\[
c_n(V^\vee)=c_n(E)-bl(V^\vee),
\]
as claimed.
\end{proof}

\begin{rem}
Note that in the above correspondence we have $L=\det(E)$.
\end{rem}

It is sometimes easy to verify when $(E,V^\vee)$ is a \emph{good pair}.

\begin{cor}\label{surfaces}
Let $X$ be a smooth projective variety \emph{of dimension $n$} such that $\mathrm{NS}(X)=\mathbb Z\cdot L$. Then there is a one-to-one correspondence
\[
\Bigl\{\, V\in \mathrm{Gr}(n+1,H^0(L)) \ \Big| \ \deg(\varphi_V)=d\,\Bigr\}
\ \longleftrightarrow\
\Bigl\{\, (E,V^\vee) \ \Big| \ 
\left(
\begin{aligned}
& E \text{ is stable reflexive of rank }n, \\
& c_1(E)=L,\ \ c_n(V^\vee)=d
\end{aligned}
\right)
\,\Bigr\}.
\]
In particular, for any $V\in \mathrm{Gr}(n+1,H^0(L))$ the associated kernel sheaf is stable; conversely, for any stable reflexive sheaf $E$ of rank $n$ and any $V^\vee\in \mathrm{Gr}(n+1,H^0(E))$, the pair $(E,V^\vee)$ is a good pair.
\label{stability}
\end{cor}

\begin{proof}
Consider $V\in \mathrm{Gr}(n+1,H^0(L))$. Since $\mathrm{NS}(X)=\mathbb{Z}\cdot L$, the evaluation map $V\otimes \mathcal O_X\to L$ is surjective in codimension $1$ (if the base locus had a divisorial component, then there would be an injective map from a non-trivial effective line bundle to $L$).

Now suppose, by way of contradiction, that the associated kernel sheaf $E$ is not stable. Then there is a destabilizing exact sequence
\[
\begin{tikzcd}
0 \arrow{r} & M \arrow{r} & E \arrow{r} & N \arrow{r} & 0,
\end{tikzcd}
\]
with $c_1(M)\ge c_1(E)=L$ numerically, hence $c_1(N)\le 0$ numerically. Since $E$ is generically globally generated, any quotient of $E$ is generically globally generated as well. A generically globally generated torsion-free sheaf with $c_1\le 0$ must be of the form $N\simeq \mathcal O_X^{\oplus k}$. Dualizing, we obtain a non-zero map $\mathcal O_X^{\oplus k}\to E^\vee$, contradicting \autoref{lemgoodpair}.

On the other hand, if we start with a stable reflexive sheaf $E$, stability implies $h^0(E^\vee)=0$. Consider any $V^\vee\in \mathrm{Gr}(n+1,H^0(E))$ and suppose, by way of contradiction, that the evaluation map $\mathrm{ev}:V^\vee \otimes \mathcal O_X\to E$ is not surjective in codimension $1$. If $\mathrm{Im}(\mathrm{ev})$ has generic rank $\le n-1$, then $c_1(\mathrm{Im}(\mathrm{ev}))>0$, hence $c_1(\mathrm{Im}(\mathrm{ev}))\ge L$ numerically since $\mathrm{NS}(X)=\mathbb Z\cdot L$. This contradicts stability of $E$. If $\mathrm{Im}(\mathrm{ev})$ has generic rank $n$, then again $c_1(\mathrm{Im}(\mathrm{ev}))\ge L$, and this implies\footnote{We have an induced non-zero map $\det(\mathrm{Im}(\mathrm{ev}))\to \det(E)=L$, hence $c_1(\mathrm{Im}(\mathrm{ev}))\le c_1(E)$.}
 $c_1(\mathrm{Im}(\mathrm{ev}))=L$ and the map $\mathrm{ev}$ is surjective in codimension $1$. Indeed, if the cokernel were supported on a divisor, then $\mathrm{Im}(\mathrm{ev})$ and $E$ would not have the same first Chern class. Hence, for any stable reflexive sheaf $E$ of rank $n$ with $c_1(E)=L$ and any $V^\vee \in \mathrm{Gr}(n+1,H^0(E))$, the pair $(E,V^\vee)$ is a good pair.
\end{proof}

In the case of surfaces, the picture is simpler. In what follows, by the intersection of a family of effective zero-cycles, we mean the effective zero-cycle of maximal degree that is a subcycle of every zero-cycle in the family.

\begin{cor}\label{degreegoodpairsurfaces}
Let $S$ be a smooth surface with $\mathrm{NS}(S)=\mathbb Z \cdot L$, and let $(E,V^\vee)$ be a good pair with $c_1(E)=L$. Then:
\begin{itemize}
\item[(i)] $E$ is a vector bundle;
\item[(ii)] for a general section $s\in V^\vee$, the cycle supported on $\mathrm{Bl}(V^\vee)$ is rigid, in the following sense: there exists a non-empty Zariski open subset $U\subset V^\vee\setminus \{0\}$ such that $A'_s$ does not depend on $s\in U$. In particular,
\[
bl(V^\vee)=\deg\Bigl(\bigcap_{s\in U} Z_{\mathrm{cycle}}(s)\Bigr)=\mathrm{deg}(A'_s);
\]
\item[(iii)] the degree of $\varphi_V$ satisfies
\[
\deg(\varphi_V)=d=c_2(E)-bl(V^\vee)>0.
\]
\end{itemize}
\end{cor}

\begin{proof}
(i) On a smooth surface, any reflexive sheaf is locally free; hence $E$ is a vector bundle.

(ii) The locus $\mathrm{Bl}(V^\vee)$ is zero-dimensional. For $s\in V^\vee$ general, the part of the zero-cycle $Z_{\mathrm{cycle}}(s)$ supported on $\mathrm{Bl}(V^\vee)$ is therefore rigid in the sense above. This implies the stated formula for $bl(V^\vee)$.

(iii) The degree computation follows from \autoref{wow}, together with the fact that any section of $E$ vanishes along a zero-dimensional subscheme: otherwise it would induce an injective map from a non-trivial effective line bundle to $L$. The inequality $d>0$ follows from the indecomposability of $L$.
\end{proof}

We end this section with a couple of remarks.

\begin{rem}\label{schemetheoreticbaselocus}
In (ii), the scheme-theoretic structure of the vanishing locus of $s\in V^\vee$ at points of $\mathrm{Bl}(V^\vee)$ may vary with $s$. Moreover, the degree of the associated cycle supported on $\mathrm{Bl}(V^\vee)$ may jump for special choices of $s\in V^\vee$.
\end{rem}

\begin{rem}\label{remark:secant-interpretation}
In the surface case, the length of the base locus satisfies $\mathrm{length}(\mathrm{Bl})=L^2-c_2(E)$. Thus $c_2(E)$ measures how \emph{secant} the linear space
\[
\mathbb P\!\bigl(\ker(H^0(L)^\vee\to V^\vee)\bigr)\subset \mathbb P H^0(L)^\vee
\]
is with respect to $S$ (scheme-theoretically). On the other hand, $bl(V^\vee)$ measures the number of points (counted with multiplicity) at which this linear space is tangent to $S$.
\end{rem}

\section{Brill--Noether loci for the very general $K3$ surface of genus $4,5,6$}\label{Section2}
In this section $(S,L)$ will denote a polarized smooth surface. We will use the framework provided by the previous section to study the loci
\[
W^2_d(S,L)=\Bigl\{\, V \in \mathrm{Gr}(3,H^0(L)) \ \Big| \ \deg(\varphi_V)\le d \Bigr\}.
\]
In the last section of the paper we will see that, as soon as $\mathrm{NS}(S)=\mathbb Z\cdot L$, these loci can always be endowed with a natural scheme structure. Under this hypothesis, points of $W^2_d(S,L)$ correspond via \autoref{wow} to pairs $(E,V^\vee)$, where $E$ is a stable rank $2$ vector bundle by \autoref{stability}. 

We begin by studying $K3$ surfaces of genus $4$, $5$, and $6$ with Picard rank one. We recall the following preliminary result, essentially due to Mukai (see \cite[Corollary~0.2]{mukai84} or \cite[Corollary~2.5]{mvbnd}). It can also be viewed as a consequence of \cite[Theorem, p.~80]{LAZ}, which shows that curves with indecomposable class on a very general $K3$ surface do not carry linear series with negative Brill--Noether number. In particular, the inequality below follows from the condition $\rho(g,r-1,d)\ge 0$. We include a proof, as we could not find a sufficiently elementary reference. For a modern treatment in the more general setting of Bridgeland-stable objects, see \cite[Theorem~6.10]{bayer}, which also implies the theorem below.

\begin{thm}[Mukai]\label{mukki}
Let $(S,L)$ be a  $K3$ surface of genus $g$ with $\mathrm{Pic}(S)=\mathbb Z\cdot L$. A stable vector bundle $F$ of rank $r$ with $c_1(F)=L$ and $c_2(F)=d$ exists if and only if
\[
d\ge \frac{(r-1)(g+r)}{r}.
\]
Moreover, if equality holds, there exists a unique vector bundle with these invariants.
\end{thm}

\begin{proof}
If such a bundle exists, it must satisfy the inequality by \cite[Corollary~2.5]{mvbnd}. In the proof of \autoref{cohomology} below we only need the \emph{non-existence} direction of the theorem; we may therefore use \autoref{cohomology} to establish the remaining implications.

For the existence direction, assume $\rho(g,r-1,d)\ge 0$, equivalently $d\ge \frac{(r-1)(g+r)}{r}$. Choose a smooth curve $C\in |L|_{\mathrm{sm}}$ which is Brill--Noether--Petri general (see \cite{LAZ} together with \cite{kleiman1972existence}) and a globally generated line bundle $A\in W^{r-1}_d(C)$. Consider the Lazarsfeld--Mukai bundle $E_{C,A}$, defined by the exact sequence
\[
\begin{tikzcd}
0 \arrow{r} & E_{C,A}^\vee \arrow{r} & H^0(A)\otimes \mathcal O_S \arrow{r} & A \arrow{r} & 0 .
\end{tikzcd}
\]
From the defining sequence we get $H^0(E_{C,A}^\vee)=0$ and $E_{C,A}$ generically globally generated, which implies that $E_{C,A}$ is stable: indeed, a destabilizing quotient $E_{C,A}\twoheadrightarrow Q$ would be generically globally generated (by $H^0(A)^\vee\subset H^0(E_{C,A})$) and satisfy $c_1(Q)\le 0$, hence $Q\simeq \mathcal O_S^{\oplus k}$. This would yield  $0\neq H^0(Q^\vee)\subset H^0(E_{C,A}^\vee)$, a contradiction. The bundle $E_{C,A}$ is the stable vector bundle with the desired invariants.

Finally, assume $\rho(g,r-1,d)=0$ and let $E,F$ be stable bundles with the same invariants $c_1^2=2g-2$, $c_2=d$, and $\mathrm{rk}=r$. Then one computes
\[
\chi(E\otimes F^\vee)=2-2\rho(g,r-1,d)=2.
\]
 In particular, there exist non-zero morphisms $\varphi:E\to F$ and $\psi:F\to E$. Since $E$ is simple, the endomorphism $\varphi\circ \psi$ is either $0$ or an isomorphism. Suppose, by way of contradiction, that $\varphi\circ \psi=0$. Then $\varphi$ is not of maximal rank. By \autoref{cohomology} (applied with $d'=d$), both $E$ and $F$ are globally generated. Since both $\coker(\varphi)$ and $\mathrm{Im}(\varphi)$ are (generically torsion-free) quotients of generically globally generated bundles, they have $c_1\ge 0$. As $L=\mathrm{det}(F)=\mathrm{det}(E)$ is indecomposable, this forces either $c_1(\coker(\varphi))=0$ or $c_1(\mathrm{Im}(\varphi))=0$. In either case one gets a contradiction, since $h^2(E)=h^2(F)=0$ and  (again, a generically globally generated sheaf with $c_1=0$ has torsion free part isomorphic to $\mathcal O_S^{\oplus k}$). Therefore, $\varphi\circ\psi$ is an isomorphism, and $E\simeq F$, proving uniqueness.
\end{proof}
\begin{rem}
One could also obtain a (somewhat more involved) proof of the \emph{non-existence} direction of \autoref{mukki} by reversing the Lazarsfeld--Mukai construction and combining it with the \emph{non-existence} part of Brill--Noether theory for curves on $K3$ surfaces (i.e. starting with a stable bundle with fixed invariants and a suitable space of sections construct a curve in $|L|_{sm}$ together with a special linear series). This is similar, in spirit, to the approach taken in the present paper to study linear systems $V\in \mathrm{Gr}(k,H^0(L))$.
\end{rem}
\autoref{mukki} has the following corollary of independent interest, which is likely known to experts.

\begin{cor}\label{cohomology}
Let $F$ be a stable vector bundle on a $K3$ surface $(S,L)$ of Picard rank one, of rank $r$, with $c_1(F)=L$ and $c_2(F)=d$. Then for any $d'< \frac{r(g+r+1)}{r+1}.$ and any $\xi \in S^{[d'-d]}$ one has
\[
h^1(F\otimes \mathcal I_\xi)=0.
\]
\end{cor}
\begin{proof}
Suppose, by way of contradiction, that there exists $\xi \in S^{[d'-d]}$ with $h^1(F\otimes \mathcal I_\xi)\ge 1$ with  $d'<\frac{r(g+r+1)}{r+1}$. Then there is a non-trivial extension
\[
\begin{tikzcd}
0 \arrow{r} & F^\vee \arrow{r} & G \arrow{r} & \mathcal I_\xi \arrow{r} & 0,
\end{tikzcd}
\]
so that $c_2(G)=c_2(F)+\deg(\xi)=d'$ and $c_1(G)=-c_1(F)$.

We claim that $G^\vee$ is stable. Indeed:
\begin{itemize}
\item replacing $G$ by its double dual, we may assume that $G$ is locally free (this amounts to replacing $\xi$ by a subscheme $\xi'\subseteq \xi$, and hence only decreases $c_2(G)$, so in particular $c_2(G)\le d'$ still holds);
\item if $K\hookrightarrow G$ is a destabilizing locally free subsheaf, then $c_1(K)\ge 0$. Moreover, the induced map $K\to \mathcal I_\xi$ cannot be injective (otherwise $K\simeq \mathcal O_S$ would split the extension);
\item thus $\ker(K\to \mathcal I_\xi)$ is a non-zero locally free subsheaf of $F^\vee$ with $c_1\ge 0$, contradicting the stability of $F^\vee$.
\end{itemize}
Therefore $G^\vee$ is stable (of rank $r+1)$, and \autoref{mukki} yields the the desired inequality, contradicting $d'<\frac{r(g+r+1)}{r+1}$.
\end{proof}
 \begin{prop}
Let $(S,L)$ be a $K3$ surface of genus $4$ with $\mathrm{Pic}(S)=\mathbb Z\cdot L$. Then $W^2_3(S,L)=\mathbb P^3.$ 
\end{prop}

\begin{proof}
By \autoref{mukki}, there is a unique stable rank $2$ vector bundle $E$ on $S$ with $c_1(E)=L$ and $c_2(E)=3$. We claim that, via the correspondence of \autoref{surfaces}, one has
\[
W^2_3(S,L)=\mathrm{Gr}(3,H^0(E))=\mathbb P H^0(E)^\vee\simeq \mathbb P^3.
\]
On the one hand, we clearly have $\mathrm{Gr}(3,H^0(E))\subset W^2_3(S,L)$ through the correspondence $(E,V^\vee)\mapsto (L,V)$. 

On the other hand, giving a degree-$3$ projection $S\subset \mathbb P^4 \dashrightarrow \mathbb P^2$ is the same as projecting away from a line $l\subset \mathbb P^4$ which is $3$-secant to $S$. The evaluation exact sequence associated with this rational map then takes the form
\[
\begin{tikzcd}
0 \arrow{r} & E_l^\vee \arrow{r} & V\otimes \mathcal O_S \arrow{r} & L \otimes \mathcal I_{l \cap S} \arrow{r} & 0.
\end{tikzcd}
\]
The kernel $E_l$ is a stable rank $2$ vector bundle with
\[
c_2(E_l)=L^2-\deg(l \cap S)=6-3=3
\qquad\text{and}\qquad
c_1(E_l)=L
\]
(see \autoref{surfaces}). Therefore, by the uniqueness part of \autoref{mukki}, $E_l\simeq E$, and hence $V^\vee\in \mathrm{Gr}(3,H^0(E))$, as desired.
\end{proof}

We now move to $K3$ surfaces of genus $6$.

\begin{thm}\label{genus6}
Let $(S,L)$ be a  $K3$ surface of genus $6$ with $\mathrm{Pic}(S)=\mathbb Z \cdot L$. Then $W^2_3(S,L)=S$. In particular, $\mathrm{irr}(S)=\mathrm{irr}_L(S)=3$.
\end{thm}

\begin{proof}
By \autoref{mukki}, on such a surface there is a unique stable rank $2$ vector bundle $E$ with invariants $c_1(E)=L$ and $c_2(E)=4$, and one has $h^0(E)=5$. Moreover, since $E$ is globally generated by \autoref{cohomology}, we have $h^0(E\otimes \mathcal I_P)=3$ for every $P\in S$. 

By the correspondence of \autoref{surfaces}, for any $P\in S$ the pair $(E,H^0(E\otimes \mathcal I_P)^\vee)$ is a good pair and induces a map of degree
\[
d=c_2(E)-bl\bigl(H^0(E\otimes \mathcal I_P)^\vee\bigr)=4-1=3.
\]
This gives an inclusion $S\hookrightarrow W^2_3(S,L)$; we view $S$ as the image of this inclusion.

Moreover, it can be shown that any pair $(L,V)\in W^2_3(S,L)$ with $\deg(\varphi_V)\le 3$ and $V\subset H^0(L)$ gives rise to a kernel bundle satisfying $c_2(E)\le 4$; see \cite[Corollary~2.3]{morettirojas}. Since $E$ is stable by \autoref{surfaces}, \autoref{mukki} implies that $c_2(E)\ge 4$. Hence $c_2(E)=4$ and $E$ must be the unique stable rank $2$ vector bundle with $c_1(E)=L$ and $c_2(E)=4$.

Finally, in view of \autoref{degreegoodpairsurfaces}, since $d=3$ we have $bl(V^\vee)=c_2(E)-d=1$, so the sections of $V^\vee$ have a base point. Therefore $V^\vee=H^0(E\otimes \mathcal I_P)^\vee$ for some $P\in S$, and thus $(L,V)$ belongs to the image of the inclusion $S\hookrightarrow W^2_3(S,L)$. We conclude that $W^2_3(S,L)=S$.

The fact that $\mathrm{irr}(S)\neq 2$ follows from the fact that a very general $K3$ surface does not admit any non-trivial dominant rational map to $\mathbb P^2$ of degree $2$ (see, for instance, \cite[Corollary~2.12]{huybrechts2016lectures}).
\end{proof}

\begin{rem}
The geometric interpretation of this projection is transparent. We have an exact sequence
\[
\begin{tikzcd}
0 \arrow{r} & E^\vee \arrow{r} & H^0(E\otimes \mathcal I_P)^\vee \otimes \mathcal O_S \arrow{r} & L \otimes \mathcal I_{\mathrm{Bl}} \arrow{r} & 0,
\end{tikzcd}
\]
where we are considering the base locus $\mathrm{Bl}$ scheme-theoretically. It is easy to see that $\deg(\mathrm{Bl})=L^2-c_2(E)=6$. Hence we are projecting from a $6$-secant linear subspace of codimension $3$. Moreover, the map
\[
H^0(E\otimes \mathcal I_P)\otimes k(P)\to E\otimes k(P)
\]
is zero, hence $\mathcal I_{\mathrm{Bl}}\otimes k(P)$ is $3$-dimensional. This implies that the codimension-$3$ linear subspace is tangent to $S$ at $P$.
\end{rem}

\, We now study projections $S \subset \mathbb{P}^5 \dashrightarrow \mathbb{P}^2$ for a very general $K3$ surface of genus $5$.

\begin{thm}
Let $(S,L)$ be a $K3$ surface of genus $5$ such that the Picard group has rank $1$. Then there are no projections $S \subset \mathbb{P}(H^0(L)^\vee) \dashrightarrow \mathbb{P}^2$ of degree $\le 3$.
\end{thm}

\begin{proof}
Let $S$ be a $K3$ surface of genus $5$ and suppose that there exists a rational map $\varphi:S\dashrightarrow \mathbb P^2$ of degree $3$ induced by a $3$-dimensional subspace $V\subset H^0(L)$. By \autoref{mukki} and \cite[Corollary~2.3]{morettirojas}, the associated kernel bundle $E$ satisfies $c_2(E)=4$. By \autoref{cohomology} we have $h^1(E)=0$, hence by Riemann--Roch $h^0(E)=\chi(E)=4$.
Again \autoref{cohomology} implies that $h^1(E\otimes \mathcal I_P)=0$ and $h^0(E\otimes \mathcal I_P)=2$ for every $P\in S$. Hence, for any $V\in \mathrm{Gr}(3,H^0(E))$ the induced map $\varphi_V$ has degree $4$ by \autoref{degreegoodpairsurfaces}.

\end{proof}

\begin{rem}
The content of the previous theorem can be rephrased as follows: a $4$-secant plane to $S\subset \mathbb P^5$ is never tangent to $S$ (more precisely, its intersection with $S$ is always a local complete intersection subscheme).
\end{rem}

We can also describe the Brill--Noether locus $W^2_4(S,L)$ in this case.

\begin{cor}\label{cor:W24}
The locus $W^2_4(S,L)\subset \mathrm{Gr}(3,H^0(L))$ consists of two disjoint components, one birational to $S$ and the other birational to a $\mathbb P^3$-bundle over $\mathcal M$, where $\mathcal M$ denotes the moduli space of stable rank $2$ vector bundles on $S$ with $c_1=L$ and $c_2=4$.
\end{cor}

\begin{proof}
The image of the map
\begin{align}
S &\hookrightarrow W^2_4(S,L)\subset \mathrm{Gr}(3,H^0(L))\label{eq:S-into-W24}\\
P &\longmapsto H^0(L\otimes \mathcal I_P^2)\nonumber
\end{align}
is the component birational to $S$. It parametrizes projections from planes in $\mathbb P^5$ that are tangent to $S$.

For the other component, let $\mathcal M$ be the moduli space of stable rank $2$ vector bundles on $S$ with $c_1=L$ and $c_2=4$, and let $\mathcal G\to \mathcal M$ be the relative Grassmannian of hyperplanes,
\[
\mathcal G_E=\mathrm{Gr}(3,H^0(E))\simeq \mathbb P^3.
\]
By the correspondence of \autoref{surfaces}, we obtain a morphism
\begin{align}
\mathcal G &\longrightarrow W^2_4(S,L)\label{eq:G-to-W24}\\
(E,V) &\longmapsto V^\vee.\nonumber
\end{align}
This component parametrizes maps whose scheme-theoretic base locus has degree $4$.

The two components are disjoint. Indeed, if they intersected we would obtain a tangent plane to $S$ meeting $S$ in at least two points. Projecting from this plane would give a rational map $S\dashrightarrow \mathbb P^2$ of degree $\le 3$, contradicting the previous theorem.

Finally, we show that there are no other components by considering the possible degrees of the base locus. If the base locus has degree $\le 2$, then $\deg(\varphi_V)\ge 6$. If the base locus has degree $3$, then it must be of the form $\mathrm{Spec}\,\mathcal O_S/\mathcal I_P^2$; otherwise the map would have degree $\ge 5$. In this case the rational map belongs to the component birational to $S$.

If the base locus has degree $4$, then the associated kernel is a stable rank $2$ vector bundle with $c_2=4$, hence the map lies in the component coming from \eqref{eq:G-to-W24}. Finally, the base locus cannot have higher degree, because there are no stable rank $2$ vector bundles with $c_2\le 3$ by \autoref{mukki}. This completes the proof.
\end{proof}

As mentioned in the introduction, the interested reader may find further computations of the loci $W^2_d(S,L)$ for $K3$ surfaces up to genus $14$ in \cite{morettirojas}.
\section{The case of the $(1,6)$-polarized abelian surface}\label{section3}

In this section we consider a very general $(1,6)$-polarized abelian surface $(A,L_A)$. Every such surface admits a presentation as a quotient
\[
(A,L_A)\simeq (B/\langle b\rangle,\,L_A),
\]
where $(B,L_B)$ is a $(1,3)$-polarized abelian surface and $b\in B[2]$ is a non-trivial $2$-torsion point.\footnote{In fact, there are three such presentations, corresponding to the three non-trivial $2$-torsion points in the kernel of the dual polarization on $\mathrm{Pic}^0(A)$.}
Starting from $(B,L_B)$ and $b\in B[2]$, let
\[
\pi:B\longrightarrow A:=B/\langle b\rangle
\]
be the quotient map. The polarization $L_A$ is induced by pushing forward divisors in $|L_B|$; for instance, by the projection formula one checks that $L_A^2=12$. Conversely, given $(A,L_A)$, choosing a non-trivial $2$-torsion point in the kernel of the dual polarization recovers such a presentation.

{\bf Assumption.} We  assume that $L_B$ contains no invertible subsheaf of the form $\mathcal O_B(nF)$ with $n\ge 2$ for some elliptic curve $F\subset B$, and that divisors in $|L_B|$ do not contain elliptic curves invariant under translation by $b$. These conditions hold for a general abelian surface.

Our goal is to construct generically finite rational maps $A\dashrightarrow \mathbb P^2$ of degree $3$. We will use the rank-$2$ vector bundle
\[
E:=\pi_*L_B
\]
to do so. We begin by recording some basic properties of $E$. If one assumes $\mathrm{NS}(A)=\mathbb Z\cdot L_A$, then $E$ is stable, and the proofs of \autoref{lem:pushforwardE} and \autoref{cor:irr16} can be simplified substantially in view of \autoref{surfaces} and \autoref{degreegoodpairsurfaces}.

\begin{lem}\label{lem:pushforwardE}
The push-forward $E=\pi_*L_B$ is a rank $2$ vector bundle with invariants $c_1(E)=L_A$ and $c_2(E)=3$. Moreover, $E$ is globally generated in codimension $1$, $h^0(E^\vee)=0$, the determinant map
\[
\bigwedge^2 H^0(E)\longrightarrow H^0(L_A)
\]
is injective and, in particular, the couple $(E,H^0(E))$ is a good pair.
\end{lem}

\begin{proof}
We observe that
\[
\pi^*E\simeq L_B\oplus t_b^*L_B,
\]
so the computation of the Chern classes is immediate. Moreover, under the identification
\[
i:H^0(L_B)\xrightarrow{\ \sim\ } H^0(E),
\]
for $s\in H^0(L_B)$ we have
\[
\pi^*i(s)=(s,t_b^*s),
\qquad\text{hence}\qquad
Z(\pi^*i(s))=Z(s)\cap t_b^*Z(s).
\]

\smallskip
\noindent\emph{Injectivity of the determinant map.}
Since $h^0(E)=3$, every element of $\bigwedge^2 H^0(E)$ is a simple tensor. If $i(s_1)\wedge i(s_2)=0$, then $i(s_1)$ and $i(s_2)$ span an invertible subsheaf $M\subset E$. In particular, $Z(i(s_1))$ moves in a positive-dimensional linear system. Pulling back to $B$, the subsheaf $\pi^*M\subset \pi^*E\simeq L_B\oplus t_b^*L_B$ yields an invertible subsheaf of $L_B$ with two independent sections, whose vanishing divisor is invariant under translation by $b$ (since it is the pullback of a divisor on $A=B/\langle b\rangle$).

This can only occur if $M^2=0$. Indeed, otherwise $M$ would descend to a line bundle on $A$ of type $(a,b)$ with $2ab\ge 4$, and its pullback cannot be an invertible subsheaf of $L_B$, which has type $(1,3)$. Therefore $\pi^*M\simeq \mathcal O_B(nF)$ for some elliptic curve $F\subset B$ and some $n\ge 2$ (since $h^0(M)\ge 2$), contradicting our assumptions on $L_B$. This proves injectivity.

\smallskip
\noindent\emph{Zero-dimensional vanishing for a general section.}
We claim that for a general $s\in H^0(L_B)$ the section $i(s)\in H^0(E)$ has vanishing locus of dimension $0$. Otherwise, for every $s\in H^0(L_B)$ the scheme $Z(\pi^*i(s))=Z(s)\cap Z(t_b^*s)$ would contain a divisor $D_s$, and hence we would obtain an injection $\mathcal O_B(D_s)\hookrightarrow L_B$. To simplify notation, we suppress the dependence on $s$ and write $D=D_s$.

Since $D$ is $t_b$-invariant, there are two possibilities:
\begin{itemize}
\item $D$ induces a $(1,2)$ polarization; or
\item $D$ is an elliptic curve (it cannot be a union of elliptic curves by our assumptions on $L_B$). In this case $D$ does not depend on $s$, since $s$ varies in $\mathbb P^2$ and $D$ has no non-trivial linearly equivalent translates.
\end{itemize}
In the first case, a general section of $H^0(\mathcal O_B(D))\subset H^0(L_B)$ is not $t_b$-invariant, hence $D$ cannot occur as the divisorial part of $Z(s)\cap Z(t_b^*s)$ for all $s$. In the second case, if every section vanished along a fixed $t_b$-invariant elliptic curve $D$, then $D$ would be a fixed component of $|L_B|$, and $H^0(L_B)=H^0(L_B(-D))$. This forces $h^1(L_B(-D))>0$ (since $L_B(-D)^2<6$), hence any divisor in $|L_B(-D)|$ is a union of elliptic curves, contradicting our assumptions. This proves the claim.

\smallskip
\noindent\emph{A degree-$3$ map and global generation in codimension $1$.}
Set
\[
E_1:=\bigl(\mathrm{Im}(H^0(E)\otimes \mathcal O_A\to E)^\vee\bigr)^\vee.
\]
Via the correspondence of \autoref{wow}, we obtain a non-degenerate rational map
\[
\varphi:A\dashrightarrow \mathbb P H^0(E)=\mathbb P H^0(E_1).
\]
Since a general section of $H^0(E)$ has zero-dimensional vanishing locus of length $c_2(E)=3$ as a section of $E$, it follows that a general section of $H^0(E_1)$ has zero-dimensional vanishing locus of length $\le 3$. Hence, by the fiber description in \autoref{wow}, the degree of $\varphi$ is $\le 3$. By the main theorem of \cite{AlzatiPirola}, the degree cannot be $1$ or $2$, so either $\varphi$ has degree $3$ or its image is a (non-degenerate) curve $C\subset \mathbb P H^0(E_1)$.

We claim that the latter does not occur. Indeed, if $\mathrm{Im}(\varphi)=C$, then pulling back a line $l\subset \mathbb P H^0(E_1)$ would produce a divisor $D_l\subset A$ which is numerically a positive multiple $k\cdot D'$, where $k=\deg(C)$ and $D'$ is the pullback of a point of $C$. If $k=2$, then $C$ is rational and $h^0(\mathcal O_A(D'))\ge 2$, so either $D'^2\ge 4$ (contradicting $\mathcal O_A(2D')\subset L_A$) or $D'^2=0$ and $D'$ is a union of (at least) two elliptic curves, contradicting our assumptions after pulling back to $B$. The argument is analogous for $k\ge 3$ when $C$ is rational.

If $C$ is not rational, then $C$ is an elliptic curve and the fibers of $A\to C$ are unions of $k'$ elliptic curves. Writing $D'=k'\cdot F$, if $k'\ge 2$ then $\pi^*\mathcal O_A(D')\subset L_B$ gives an invertible subsheaf of $L_B$, contradicting our assumptions. If $k'=1$, then by the fiber description in \autoref{wow} a section $[i(s)]\in C$ vanishes along an elliptic curve as a section of $E_1$, hence also as a section of $E$. Pulling back to $B$ and using $\pi^*Z(i(s))=Z(s)\cap Z(t_b^*s)$, we conclude that $s$ either vanishes along a $t_b$-invariant elliptic curve, or vanishes along two elliptic curves exchanged by translation by $b$, again contradicting our assumptions on $L_B$. This proves the claim.

Therefore the map associated to $(E_1,H^0(E))$ has degree $3$. To conclude the proof, we show that $E_1=E$, i.e.\ that $E$ is globally generated in codimension $1$. The group $K(L_B)\simeq \mathbb Z_3\times \mathbb Z_3$ acts on $H^0(L_B)$, hence on $H^0(E)=H^0(E_1)$ and therefore on $\bigwedge^2 H^0(E_1)^\vee\subset H^0(\det(E_1))$ (see also \cite[Theorem~2.2]{grossi2025explicitclasslagrangiansurfaces}). Since $E_1$ is globally generated in codimension $1$, the base locus $\mathrm{Bl}$ of $|\bigwedge^2H^0(E_1)|$ has codimension $2$ and is invariant under $K(L_B)$. As $K(L_B)$ acts on $A$ by translations (hence has no fixed points), we deduce that $|\mathrm{Bl}|\ge 9$. Since $\deg(\varphi)\ge \det(E_1)^2-|\mathrm{Bl}|$, we obtain $\det(E_1)^2\ge 12$.

Finally, since $E_1\subset E$, we have an inclusion $\det(E_1)\subset \det(E)=L_A$. On an abelian surface, if $M\subset N$ are effective line bundles with $M^2\ge N^2>0$, then $M=N$. Hence $\det(E_1)=L_A$, and since $E_1\subset E$ with both locally free of the same rank, we conclude that $E_1=E$.
\end{proof}

\begin{cor}\label{cor:irr16}
The degree of irrationality of a general $(1,6)$-polarized abelian surface is $3$.
\end{cor}

\begin{proof}
Since $(E,H^0(E))$ is a good pair, the correspondence of \autoref{wow} yields a rational map $f:A\dashrightarrow \mathbb P^2$ of degree $\le 3$. We claim that $\deg(f)>0$. Indeed, if the image were a (non-degenerate) curve $C\subset \mathbb P^2$, then we would obtain
\[
L_A\simeq f^*(z)
\]
for some $z\in \mathrm{Pic}(C)$ of degree at least $2$, contradicting the numerical primitivity of $L_A$. Finally, the degree cannot be $2$ by the main theorem of \cite{AlzatiPirola}. Hence $\deg(f)=3$.
\end{proof}
\begin{rem}\label{abeliansurface}
In \cite{grossi2025explicitclasslagrangiansurfaces} the geometry of the degree $3$ rational map $A\dashrightarrow \mathbb P^2$ is studied in detail, and the invariants of its Galois closure (a Lagrangian surface in $A\times A$) are computed. A similar program was carried out for $(1,2)$-polarized abelian surfaces and the degree $3$ map constructed by Tokunaga and Yoshihara (see \cite{35}) in \cite{lagrangian}.

With some extra effort it is possible to prove a theorem analogous to \autoref{456} for a very general $(1,6)$-polarized abelian surface $A\subset \mathbb P^5=\mathbb P H^0(L_A)^\vee$. One uses the geometry of stable vector bundles of rank $2$ on $A$ with $c_1=L_A$ and $c_2=3$, which turn out to be pushforwards of line bundles on $(1,3)$-polarized abelian surfaces under a degree-$2$ isogeny (see \cite[Theorem~5.8]{mu}). This implies several equivariance  properties of the degree $3$ rational map $A\dashrightarrow \mathbb P^2$ (see \cite[Section~2]{grossi2025explicitclasslagrangiansurfaces}). This yields, for such an embedding $A\subset \mathbb P^5$, the existence of exactly four planes in $\mathbb P^5$ which are $9$-secant to $A$, together with a description of their intersection with $A$.
\end{rem}

\section{The upper bound}\label{section4}
We use the construction of the previous sections to produce generically finite rational maps of low degree from very general $K3$ surfaces of genus $g$ to $\mathbb P^2$. We establish our bounds for very general $K3$ surfaces, and general specialization results (see \cite[Corollary~C]{7}) imply the same bounds for arbitrary $K3$ surfaces. Our construction improves the best previously known bounds due to Stapleton by a factor between $3$ and $6$, depending on the genus; see \cite{St}.

The following table collects some of the resulting upper bounds for small values of $g$:
\[
\begin{tabular}{l|l|l|l|l|l|l|l|l|l|l|l|l|l|l}
$g$ & 6 & 8 & 10 & 12 & 14 & 16 & 18 & 20 & 22 & 24 & 26 & 42 & 62 & $2+2n(n+1)$ \\
\hline
$\mathrm{irr}(S)\le$ & 3 & 4 & 4 & 5 & 4 & 5 & 5 & 6 & 6 & 7 & 5 & 6 & 7 & $2+n$
\end{tabular}
\]

We first deal with the case $g=2+2n(n+1)$. In view of \autoref{mukki}, the minimal possible second Chern class for a rank $2$ bundle with primitive $c_1$ is
\[
c_2(E)=2+n(n+1).
\]
We will prove that
\[
\mathrm{irr}(S)\le 2+\frac{n(n+1)}{2}-\frac{(n-1)n}{2}=2+n \sim \frac{1}{\sqrt{2}}\sqrt{g}.
\]
For this series of genera, the resulting estimate is asymptotically a $1/6$th of Stapleton's bound.
	The theorem is an easy consequence of the construction in the previous section.
	\begin{thm}
		Let $(S,L)$ be a $K3$ surface of genus $g=2+2n(n+1)$, then $\mathrm{irr}(S) \le 2+n$.
	\end{thm}
	\begin{proof}
		Consider the rank $2$ vector bundle $E$ with invariants \[
		c_1(E)=L, \quad c_2(E)=2+n(n+1),\quad  h^0(E)=3+{n(n+1)}.\]
		Since $\mathrm{colength}(\mathcal I_P^n)=\frac{n(n+1)}{2}$ for any $P\in S$ there exists a vector space $V_P^\vee\subset H^0(E \otimes \mathcal I_P^n)$ of dimension $3$.  Any section of $V_P^\vee$ vanishes at $P$ with order $n^2$ (since it is locally defined by two functions vanishing with order $n$ at $P$). We deduce that ${bl}(V_P^\vee)=\mathrm{deg}(\cap_{s\in V_P^\vee}Z_{cycle}(s))\ge n^2.$
		By   \autoref{wow} and \autoref{degreegoodpairsurfaces}
		\[
		\mathrm{deg}(\varphi_{V_P})=c_2(V_P)^\vee\le c_2(E)-n^2= 2+n,
		\]
		where $\varphi_{V_P}$ is the map $S \dashrightarrow \mathbb{P}^2$ induced by $V_P\subset H^0(L)$.
	\end{proof}
	The general case is an easy consequence. 
	\begin{cor}
		Let $(S,L)$ be a $K3$ surface of genus $g=2+2n(n+1)+k<2+2(n+1)(n+2)$ (i.e. $k<4n+4$), then \[
		\mathrm{irr}(S)\le 2+n +\ceil*{\frac k 2}- \floor*{ \frac k 4 }.
		\]
	\end{cor}
	\begin{proof}
		In view of \autoref{mukki} on $S$ there is a stable rank $2$ vector bundle $E$ with \[
		c_1(E)=L, \quad c_2(E)= 2+n(n+1)+\ceil*{\frac k 2},\quad  h^0(E)= 3+{n(n+1)}+\ceil*{\frac k 2}.\]

        Consider distinct points $P,Q_1,\dots, Q_{\floor{\frac k 4 }}\in S$ and
        
        \[
        V^\vee \in \mathrm{Gr}(3,H^0(E \otimes \mathcal I_P^n\otimes \mathcal I_{Q_1}\otimes \dots \otimes \mathcal I_{Q_{\floor{\frac k 4 }}}),\]
		then, in view of \autoref{wow} and \autoref{degreegoodpairsurfaces}, $\varphi_V$ is of degree
		\[
		\mathrm{deg}(\varphi_V)=c_2(V^\vee)\le c_2(E)-n^2-\floor*{\frac k 4}= 2+n +\ceil*{\frac k 2}- \floor*{ \frac k 4 }.\]
	\end{proof}

  \section{An application to higher-dimensional varieties}\label{section5}

In this section we provide a few applications of \autoref{wow} to higher-dimensional varieties. We begin with a sufficient criterion ensuring that a good pair $(E,V^\vee)$ induces a \emph{generically finite} rational map. Note that, by \autoref{wow}, the associated map is generically finite if and only if $c_n(V^\vee)>0$.

\begin{lem}\label{positivedegree}
Let $(E,V^\vee)$ be a good pair on an $n$-dimensional smooth projective variety $X$. Suppose that for every $k\ge 0$ the locus
\[
\Bigl\{\, s\in V^\vee \ \Big| \ Z(s)\ \text{has an irreducible component of dimension }\ge k \Bigr\}
\]
has dimension $<n-k$. Then $c_n(V^\vee)>0$, i.e.\ $\varphi_V$ induces a generically finite rational map.

In particular, if $E$ is a vector bundle and every section in $V^\vee$ has vanishing locus of dimension $0$ on $X\setminus \mathrm{Bl}(V^\vee)$, then $c_n(V^\vee)>0$.
\end{lem}

\begin{proof}
By \autoref{wow}, the fibers of $\varphi_V:X\setminus \mathrm{Bl}(V^\vee)\to \mathbb P V^\vee$ are precisely the vanishing loci of sections in $V^\vee$ restricted to $X\setminus \mathrm{Bl}(V^\vee)$. Hence the locus
\[
\bigcup_{s\in V^\vee}\bigl(Z(s)\cap (X\setminus \mathrm{Bl}(V^\vee))\bigr)
\]
covers $X\setminus \mathrm{Bl}(V^\vee)$. The conclusion follows by a dimension count using the hypothesis on the loci of sections with positive-dimensional components.
\end{proof}

\begin{prop}\label{cor:highdim-estimates}
The following estimates hold:
\begin{itemize}
\item if $(S,L)$ is a very general $K3$ surface of genus $6$, then $\mathrm{irr}(S^{[2]})\le 6$;
\item if $(A,L)$ is a very general abelian threefold of type $(1,3,12)$ or $(1,6,6)$, then $\mathrm{irr}_L(A)\le 8$;
\item if $(S,L)$ is a very general $K3$ surface of genus $10$, then $\mathrm{irr}(S^{[3]})\le 20$.
\end{itemize}
\end{prop}

\begin{proof}
\begin{itemize}
\item \emph{The case of $S^{[2]}$.}
Let $E$ be the stable rank $2$ vector bundle on $S$ with $c_1(E)=L$, $c_2(E)=4$, and $h^0(E)=5$. It canonically induces a rank $4$ vector bundle $\mathcal E$ on $S^{[2]}$, defined by
\[
\mathcal E\otimes k(\xi)\;=\; H^0(E\otimes \mathcal O_\xi)
\qquad\text{for }\xi\in S^{[2]}.
\]

Moreover, one has $c_4(\mathcal E)=6$, and there is a natural identification
\[
i:H^0(E)\xrightarrow{\ \sim\ } H^0(\mathcal E).
\]

We briefly describe the vanishing loci of sections of $\mathcal E$. If $Z(s)=\zeta\in S^{[4]}$, then
\[
Z(i(s))\subset \Bigl\{\, \xi\in S^{[2]} \ \Big| \ \xi\subset \zeta \Bigr\}.
\]
This set is finite  if $\zeta$ is curvilinear and of cardinality exactly $6$ if $\zeta$ is reduced. By \autoref{propertiesofE}, for general $s\in H^0(E)$ the subscheme $Z(s)$ is reduced (hence curvilinear), and therefore for general $s\in H^0(E)$ the vanishing locus $Z(i(s))$ is zero-dimensional of length $6$. It follows that $H^0(\mathcal E)$ generically generates $\mathcal E$ (otherwise a general section would factor through a rank $3$ subsheaf and would have vanishing locus of dimension at least $1$). In particular, the pair $\bigl((\mathrm{Im}(\mathrm{ev})^\vee)^\vee,\, H^0(\mathcal E)\bigr)$ is a good pair.

If $\zeta$ is not curvilinear, then locally $\mathcal I_\zeta=(x^2,y^2)$ at its support, and the locus
\[
\Bigl\{\, \xi\in S^{[2]} \ \Big| \ \xi\subset \zeta \Bigr\}
\]
is $1$-dimensional (it consists of length-$2$ subschemes supported at the same point). By \autoref{propertiesofE}, for every $P\in S$ the vector space $H^0(E\otimes \mathcal I_P^2)$ has dimension at most $1$, so the locus of sections vanishing along a positive dimensional subscheme is of dimension $\le 2$. Therefore \autoref{positivedegree} yields $c_4\bigl(H^0(\mathcal E)\bigr)>0$. On the other hand, since a general section vanishes at $6$ points outside the base locus, we have $c_4\bigl(H^0(\mathcal E)\bigr)\le 6$. Hence, by \autoref{wow}, the good pair $\bigl((\mathrm{Im}(\mathrm{ev})^\vee)^\vee,\, H^0(\mathcal E)\bigr)$ induces a generically finite rational map of degree $\le 6$, proving $\mathrm{irr}(S^{[2]})\le 6$.

\item \emph{The case of abelian threefolds.}
A very general abelian threefold $(A,L_A)$ of type $(1,3,12)$ (resp.\ $(1,6,6)$) admits a description as $(B,L_B)/\langle b\rangle$, where $(B,L_B)$ is a very general $(1,1,4)$ (resp.\ $(1,2,2)$) abelian threefold and $b$ is a $3$-torsion point. Let $\pi:B\to A$ be the quotient map and set $E:=\pi_*L_B$. Then $E$ is a rank $3$ vector bundle with
\[
c_1(E)=L_A,
\qquad
c_3(E)=8,
\]
since there is an isomorphism $\pi^*E\simeq L_B\oplus t_b^*L_B\oplus t_{2b}^*L_B$.
If $\mathrm{NS}(A)=\mathbb Z\cdot L_A$ (as for a very general abelian threefold), then $(E,H^0(E))$ is a good pair in view of \autoref{surfaces}. Moreover, for a very general $A$, one checks that every section $s\in H^0(E)$ has vanishing locus of dimension $0$. For instance, one can argue by degeneration. Degenerate the abelian threefold $B$ to a product $E\times S$, where $E$ is an elliptic curve and $S$ is a very general $(1,4)$-polarized abelian surface. In this situation, one checks that a section vanishes either along a zero-dimensional subscheme or along a disjoint union of smooth elliptic curves. In particular, such vanishing loci cannot occur as flat limits of vanishing loci containing a curve on a very general (simple) abelian threefold.  Therefore \autoref{positivedegree} and \autoref{wow} yield a generically finite rational map of degree $\le c_3(E)=8$.

\item \emph{The case of $S^{[3]}$.}
The argument is analogous to the case of $S^{[2]}$, using the corresponding tautological bundle on $S^{[3]}$, and we omit the details.
\end{itemize}
\end{proof}
Here is a proof of a few technical results used in the above proof.
\begin{lem}\label{propertiesofE}
Let $(S,L)$ be a $K3$ surface of genus $6$ with $\mathrm{Pic}(S)=\mathbb Z \cdot L$, and let $E$ be the stable rank $2$ vector bundle with $c_1(E)=L$ and $c_2(E)=4$. Then:
\begin{itemize}
\item for general $s\in H^0(E)$ the zero locus $Z(s)$ is reduced (hence curvilinear);
\item for any $P\in S$ one has $h^0(E\otimes \mathcal I_P^2)\le 1$.
\end{itemize}
\end{lem}

\begin{proof}
Recall that, in view of \autoref{mukki}, such a vector bundle is unique. For the first point, take a general smooth curve $C\in |L|$. By \cite{LAZ}, $C$ is Brill--Noether general, hence it carries a base-point-free $g^1_4$; let $A\in W^1_4(C)$ be the associated line bundle. Consider the Lazarsfeld--Mukai bundle
\[
E_{C,A}^\vee:=\ker\bigl(H^0(A)\otimes \mathcal O_S\to A\bigr).
\]
One checks that $h^0(E_{C,A}^\vee)=0$, which implies stability: any destabilizing subsheaf would have $c_1\ge 0$ (by $\mathrm{Pic}(S)=\mathbb Z\cdot L$) and would force $h^0(E_{C,A}^\vee)>0$. Moreover $c_1(E_{C,A})=L$ and $c_2(E_{C,A})=4$, hence $E_{C,A}\simeq E$ by uniqueness. Finally, the sections in $H^0(A)^\vee\subset H^0(E)$ vanish along a divisor in $|A|$, and a general divisor in $|A|$ is reduced. This proves the first claim.

For the second statement, suppose, by way of contradiction, that $h^0(E\otimes \mathcal I_P^2)=2$ for some $P\in S$. Then  two independent sections $s_1,s_2\in H^0(E\otimes \mathcal I_P^2)$ vanish with multiplicity at least $4$ at $P$ (their zero locus is locally given by two equations vanishing with order $2$ at $P$). Since $c_2(E)=4$ and no section of $E$ vanishes along a divisor (otherwise $E$ would be destabilized by our assumptions on $\mathrm{Pic}(S)$), these sections cannot vanish on $S\setminus P$. On the other hand, the wedge $s_1\wedge s_2$ defines a non-zero section of $\det(E)\simeq L$ vanishing along a curve $C\in |L|$. Pick $Q\in C\setminus\{P\}$. Since $s_1$ and $s_2$ are linearly dependent at $Q$, a suitable linear combination vanishes at $Q$ and with multiplicity $\ge 4$ at $P$, contradicting $c_2(E)=4$.
\end{proof}

\begin{rem}
The cases of abelian varieties and Hilbert schemes of $K3$ surfaces can be generalized to higher dimensions. In the abelian threefold case, the kernel bundle is semi-homogeneous, i.e.\ a pushforward of a line bundle under an isogeny (as for the $(1,6)$-polarized abelian surface; see \autoref{abeliansurface}). This implies, in particular, equivariance and symmetry properties for the rational maps we construct. For more details on semi-homogeneous vector bundles on abelian varieties, see the work of Mukai \cite{mu}. 
\end{rem}

As the proof of \autoref{cor:highdim-estimates} illustrates, it can be nontrivial to use \autoref{wow} to compute the degree of the rational map associated with a good pair $(E,V^\vee)$. In \cite{Morettirojas2} it is shown that, if $E$ is a globally generated vector bundle on a smooth projective variety $X$ of dimension $n$, then for general $V^\vee\in \mathrm{Gr}(n+1,H^0(E))$ the rational map associated with $(E,V^\vee)$ via \autoref{wow} has degree exactly $c_n(E)$.

\section{Towards a Brill--Noether theory for surfaces}\label{section6}

We fix a smooth projective surface $S$ with $\mathrm{NS}(S)=\mathbb Z\cdot L$ and $h^0(L)\ge 3$.
The goal of this section is to show that the locus of special projections
\[
W^r_d(S,L)=\Bigl\{\,V\in \mathrm{Gr}(r+1,H^0(L)) \ \Big| \ \deg(\varphi_V)\le d \Bigr\}
\]
is naturally a scheme.

We also remark that Mendes-Lopes, Pardini and Pirola introduced Brill--Noether loci for irregular varieties in \cite{bn} as cohomology jump loci in $\mathrm{Pic}^0(S)$. Our definition is of a different nature.

\subsection{The correspondence}

Let $V\subset H^0(L)$ be an $(r+1)$-dimensional subspace inducing a rational map
\[
\varphi_V:S\dashrightarrow \mathbb P^r=\mathbb P(V^\vee),
\]
and let $E^\vee:=\ker\bigl(V\otimes \mathcal O_S\to L\bigr)$ be the associated kernel bundle, so that we have an exact sequence
\[
\begin{tikzcd}
0 \arrow{r} & E^\vee \arrow{r} & V\otimes \mathcal O_S \arrow{r} & L.
\end{tikzcd}
\]
(Recall that on a surface, reflexive implies locally free.)

\begin{defn}\label{def:goodrpair}
A pair $(E,V^\vee)$, where $E$ is a rank $r$ vector bundle and $V^\vee\in \mathrm{Gr}(r+1,H^0(E))/\mathrm{Aut}(E)$, is a \emph{good $r$-pair} if $V^\vee$ generates $E$ in codimension $1$ and $H^0(E^\vee)=0$.
\end{defn}

We have the following analogue of \autoref{surfaces}.

\begin{lem}\label{goodd}
Let $S$ be a smooth surface with $\mathrm{NS}(S)=\mathbb Z\cdot L$. Then there is a one-to-one correspondence
\[
\Bigl\{\, V\in \mathrm{Gr}(r+1,H^0(L)) \Bigr\}
\ \longleftrightarrow\
\Bigl\{\, (E,V^\vee)\ \Big| \ E \text{ stable of rank }r,\ c_1(E)=L,\ V^\vee\in \mathrm{Gr}(r+1,H^0(E)) \Bigr\}.
\]
\end{lem}

\begin{proof}
Same proof as \autoref{surfaces}.
\end{proof}

Let $(E,V^\vee)$ be a good $r$-pair and let $W\in \mathrm{Gr}(r-1,V^\vee)$. Consider the evaluation map
\[
\mathrm{ev}_W:W\otimes \mathcal O_S\longrightarrow E.
\]
We claim that $\mathrm{ev}_W$ is injective. 
If, by contradiction, $\mathrm{ev}_W$ is not injective we would have $c_1(\mathrm{Im}(\mathrm{ev}_W))\ge L\in \mathrm{NS}(S)$, so that $\mathrm{Im}(\mathrm{ev}_W)$ would destabilize the (stable) vector bundle $E$. Hence $\mathrm{ev}_W$ is injective and we obtain an exact sequence
\[
\begin{tikzcd}
0\arrow{r} & W\otimes \mathcal O_S \arrow{r} & E \arrow{r} & L\otimes \mathcal I_\xi \arrow{r} & 0,
\end{tikzcd}
\]
for some $\xi \in S^{[c_2(E)]}$.

\begin{defn}\label{dff}
Let $(E,V^\vee)$ be a good $r$-pair with $c_1(E)=L$. We define:
\begin{itemize}
\item the base locus $\mathrm{Bl}(V^\vee)\in \mathrm{Sym}^{L^2-c_2(E)}(S)$ as the (zero-dimensional) locus where $V^\vee$ does not generate $E$;
\item for $W\in \mathrm{Gr}(r-1,V^\vee)$,
\begin{align*}
Z(W) &:= \mathrm{Spec}\,\mathcal O_S/\mathrm{coker}(\mathrm{ev}_W)\ \in S^{[c_2(E)]},\\
Z_{\mathrm{cycle}}(W) &:= \mathrm{HC}(Z(W))\ \in \mathrm{Sym}^{c_2(E)}(S),
\end{align*}
where $\mathrm{HC}:S^{[c_2(E)]}\to \mathrm{Sym}^{c_2(E)}(S)$ is the Hilbert--Chow morphism;
\item $c_2(V^\vee)$ as the degree of the part of $Z_{\mathrm{cycle}}(W)$ supported on $S\setminus \mathrm{Bl}(V^\vee)$ for general $W\in \mathrm{Gr}(r-1,V^\vee)$;
\item $bl(V^\vee):=\deg\Bigl(\bigcap_{W\in \mathrm{Gr}(r-1,V^\vee)} Z_{\mathrm{cycle}}(W)\Bigr)$.
\end{itemize}
\end{defn}

\begin{rem}\label{rigidity}
Since $\mathrm{Bl}(V^\vee)$ is zero-dimensional, we have
\[
bl(V^\vee)=\deg\bigl(Z_{\mathrm{cycle}}(W)\cap \mathrm{Bl}(V^\vee)\bigr)
\]
for general $W\in \mathrm{Gr}(r-1,V^\vee)$. We also remark that $Z(W)$ is well defined for every $W$ because $\mathrm{ev}_W$ is injective for every $W\in \mathrm{Gr}(r-1,V^\vee)$.
\end{rem}

The above constructions globalize to the relative setting. More precisely, let:
\begin{itemize}
\item $\mathcal M_{r,c}$ be the moduli space of stable vector bundles of rank $r$ with $c_1=L$ and $c_2=c$;
\item $\mathcal G_{r,c}\to \mathcal M_{r,c}$ be the relative Grassmannian whose fiber over $E\in \mathcal M_{r,c}$ is $\mathrm{Gr}(r+1,H^0(E))$;\footnote{The construction can be made locally analytically on $\mathcal M_{r,c}$, since universal bundles exist locally analytically. It can also be globalized even when no global universal bundle exists; for the $K3$ case see \cite{mvbnd}.}
\item $\mathcal W_{r,c}\to \mathcal G_{r,c}$ be the relative Grassmannian whose fiber over $(E,V^\vee)$ is $\mathrm{Gr}(r-1,V^\vee)$;
\item $\mathcal G_r:=\bigcup_{c}\mathcal G_{r,c}$ and $\mathcal W_r:=\bigcup_{c}\mathcal W_{r,c}$.
\end{itemize}
By \autoref{goodd}, the space $\mathcal G_r$ parametrizes good $r$-pairs and comes equipped with a natural bijection
\[
\mathcal G_r \longrightarrow \mathrm{Gr}(r+1,H^0(L)).
\]
Moreover, \autoref{dff} globalizes to give morphisms
\[
Z:\mathcal W_{r,c}\to S^{[c]},\qquad
Z_{\mathrm{cycle}}:\mathcal W_{r,c}\to \mathrm{Sym}^{c}(S),\qquad
\mathrm{Bl}:\mathcal G_{r,c}\to \mathrm{Sym}^{L^2-c}(S).
\]

\smallskip
The construction is based on the following observation, which generalizes \autoref{wow} in the surface case.

\begin{cor}\label{c!}
Let $\varphi_V:S \dashrightarrow \mathbb{P}^r=\mathbb P(V^\vee)$ be a rational map with kernel bundle $E$ and base locus $\mathrm{Bl}(V^\vee)$. Then for any $W\in \mathrm{Gr}(r-1,V^\vee)$ one has an identification
\[
\varphi_V^{-1}\bigl(\mathbb P(W)\bigr)\cap (S\setminus \mathrm{Bl}(V^\vee))
=
Z(W)\cap (S\setminus \mathrm{Bl}(V^\vee)).
\]
In particular,
\[
\deg(\varphi_V)=c_2(V^\vee)=c_2(E)-bl(V^\vee).
\]
\end{cor}

\begin{proof}
We abuse notation and denote by $\varphi_V$ the induced morphism $\varphi_V:S\setminus \mathrm{Bl}(V^\vee)\to \mathbb P(V^\vee)$.
The Euler sequence
\[
\begin{tikzcd}
0\arrow{r} & \mathcal O_{\mathbb P(V^\vee)}(-1) \arrow{r} & V^\vee \otimes \mathcal O_{\mathbb P(V^\vee)} \arrow{r} & T_{\mathbb P(V^\vee)}(-1) \arrow{r} & 0
\end{tikzcd}
\]
pulls back to
\[
\begin{tikzcd}
0 \arrow{r} & L^\vee_{|_{S\setminus \mathrm{Bl}(V^\vee)}} \arrow{r} &
V^\vee \otimes \mathcal O_{S\setminus \mathrm{Bl}(V^\vee)} \arrow{r} &
E_{|_{S\setminus \mathrm{Bl}(V^\vee)}} \arrow{r} & 0,
\end{tikzcd}
\]
so that $E_{|_{S\setminus \mathrm{Bl}(V^\vee)}}\simeq \varphi_V^*T_{\mathbb P(V^\vee)}(-1)$ and pullback induces an identification,
\[
\varphi_V^*:H^0\!\left(T_{\mathbb P(V^\vee)}(-1)\right)\xrightarrow{\ \sim\ } V^\vee \subset H^0(E).
\]
and, more generally, an inclusion
\[
\mathrm{Gr}\bigl(r-1, H^0(T\mathbb P(V^\vee)(-1))\bigr)\subset \mathrm{Gr}\bigl(r-1,H^0(E)\bigr).
\]
For a general $W\in \mathrm{Gr}\bigl(r-1,H^0(T_{\mathbb P(V^\vee)}(-1))\bigr)$, the map
\[
W\otimes \mathcal O_{\mathbb P(V^\vee)} \longrightarrow T_{\mathbb P(V^\vee)}(-1)
\]
drops rank along the linear subspace $\mathbb P(W)\subset \mathbb P(V^\vee)$ (and at no other point, since $c_r\!\left(T_{\mathbb P(V^\vee)}(-1)\right)=\mathcal O(1)^{r}$). Hence, for $W\in \mathrm{Gr}(r-1,V^\vee)$, the fiber of
\[
\varphi_V: S\setminus \mathrm{Bl} \longrightarrow \mathbb P(V^\vee)
\]
above $\mathbb P(W)$ is $Z_{}(W)$. The degree statement follows because $c_2(V^\vee)$ is precisely the degree of $Z_{\mathrm{cycle}}(W)$ supported outside the base locus $\mathrm{Bl}$ (see also \autoref{rigidity}), and because $\deg(\varphi_V)$ equals the degree of the pullback of a general codimension-$2$ linear subspace $\mathbb P(W)\subset \mathbb P(V^\vee)$, with $W\in \mathrm{Gr}(r-1,V^\vee)$.
\end{proof}

We are now ready to show that the loci $W^r_d(S,L)$ are quasi-projective algebraic varieties.

\begin{thm}\label{scheme}
Let $(S,L)$ be a polarized surface with $\mathrm{NS}(S)=\mathbb Z\cdot L$. Then the locus
\[
W^r_d(S,L)=\Bigl\{\, V \in \mathrm{Gr}(r+1,H^0(L)) \ \Big| \ \deg(\varphi_V)\le d \Bigr\}
\]
is a quasi-projective algebraic variety.
\end{thm}

\begin{proof}
By \autoref{goodd} we have a bijection $\mathcal G_r \to \mathrm{Gr}(r+1,H^0(L))$. We construct the loci inside $\mathcal G_r=\bigcup_c \mathcal G_{r,c}$, which filters $\mathrm{Gr}(r+1,H^0(L))$ according to the second Chern class of the kernel bundle. Let $\pi:\mathcal W_{r,c}\to \mathcal G_{r,c}$ be the relative Grassmannian. To a point $(E,V^\vee,W)\in \mathcal W_{r,c}$ we associate the two cycles $\mathrm{Bl}(V^\vee)$ and $Z_{\mathrm{cycle}}(W)$.

In view of \autoref{c!}, the point $(E,V^\vee)$ gives a rational map of degree $\le d$ if and only if for every $W\in \mathrm{Gr}(r-1,V^\vee)$ the cycle $Z_{\mathrm{cycle}}(W)$ contains at least $c-d$ points (counted with multiplicity) supported on $\mathrm{Bl}(V^\vee)$. Define
\[
\Delta_{c,d}:=\Bigl\{\, (P,Q)\in \mathrm{Sym}^c(S)\times \mathrm{Sym}^{L^2-c}(S)\ \Big|\ 
\deg(P\cap \mathrm{Supp}(Q))\ge c-d \Bigr\}.
\]
Then the locus
\[
W^r_d(S,L)_c
:=
\Bigl\{\, (E,V^\vee)\in \mathcal G_{r,c} \ \Big|\ 
\pi^{-1}(E,V^\vee)\subset (Z_{\mathrm{cycle}},\,\mathrm{Bl}\circ \pi)^{-1}(\Delta_{c,d})
\Bigr\}
\]
parametrizes good $r$-pairs with $c_2(E)=c$ and $\deg(\varphi_V)\le d$. The loci $W^r_d(S,L)_c$ are quasi-projective (by general facts on flag varieties), and therefore the finite union
\[
W^r_d(S,L)=\bigcup_c W^r_d(S,L)_c
\]
is a quasi-projective algebraic variety.
\end{proof}

	\bibliographystyle{plain}
\bibliography{refer}
\end{document}